 \documentclass[12pt]{amsart}
\usepackage{amssymb}
\usepackage{epsf}
\usepackage{xy}

\numberwithin{equation}{section}

\theoremstyle{plain}
 \newtheorem{theo}{Theorem}[section]

 \newtheorem{coro}[theo]{Corollary}

\theoremstyle{definition}
 \newtheorem{defi}[theo]{Definition}
  
 \newtheorem{exa}[theo]{Example}   
 \newtheorem{rems}[theo]{Remarks}

\newcommand{\parti}{\vdash}              
           
 \newcommand{\calA}{\mathcal{A}}
 \newcommand{\calB}{\mathcal{B}}
 
 \newcommand{\calF}{\mathcal{F}}

 \newcommand{\calT}{\mathcal{T}}

\def\Lie{{ \mathcal {L}ie}}   
\def\sb{{\hbox{sb}}}

\xyoption{all}
\CompileMatrices

\begin{document}

 \title[Basis for the free Lie algebra of the labelled rooted  trees]{A combinatorial basis for the free Lie algebra of the labelled rooted  trees}

 \author[N.~Bergeron and M.~Livernet]{Nantel Bergeron
 and Muriel Livernet}

 \address[Nantel Bergeron]
 {Department of Mathematics and Statistics\\York University\\
 Toronto, Ontario M3J 1P3\\
 CANADA}
 \email{bergeron@mathstat.yorku.ca}
 \urladdr[Nantel Bergeron]{http://www.math.yorku.ca/bergeron}

 \address[Muriel Livernet]
 {LAGA, Institute Galil\'ee \\ Universit\'e Paris13\\ Avenue
   Jean-Baptiste Cl\'ement\\ 93430 Villetaneuse\\France}
\email{livernet@math.univ-paris13.fr}
\urladdr[Muriel Livernet]{http://www.math.univ-paris13.fr/~livernet}

 \date{}

 \thanks{Bergeron supported in part by CRC and NSERC}
 

 \begin{abstract}
The pre-Lie operad can be realized as a space $\calT$ of labelled rooted trees.
 A result of F. Chapoton shows that the pre-Lie operad is a free twisted Lie algebra. That is $\calT=\Lie \circ \calF$ for some $\mathbb S$-module $\calF$. In the context of species, we construct  an explicit basis of $\calF$. This allows us to give a new proof of Chapoton's results. Moreover it permits us to show that $\calF$ forms a sub nonsymmetric operad of the pre-Lie operad $\calT$.
 \end{abstract}

 \maketitle

\vspace*{-0.3cm}\begin{figure}[!h]
\parbox{350pt}{\tiny\tableofcontents}
\end{figure}

 
 \section*{Introduction}
 
 One of the most fascinating result in the theory of operads is the
 Koszul duality between the 
Lie operad $\Lie$ and the commutative and associative operad ${\mathcal C}om$ \cite{GK94}.
   This has inspired many researcher to study this pair of operads and its refinements. 
One particular instance of this is the study of pre-Lie algebras. 
That is vector space $L$ together with a product $*$ satisfying the relation 
$$(x* y)* z-x*(y* z)=(x* z)* y-x*(z* y), \ \forall x,y,z\in L.$$
In \cite{CL01}, we have a realization of the pre-Lie operad as the space $\calT$ of labelled rooted trees.
This operad sits naturally between $\Lie$ and ${\mathcal A}s$  as we have the injective morphisms
$$\Lie\rightarrow \mathcal T \rightarrow {\mathcal A}s$$
which factor the usual injective morphism  from $\Lie$ to ${\mathcal A}s$. 

 At the level of algebras, the free pre-Lie algebra generated by a
vector space of dimension 1, is the space of unlabelled rooted
trees. Indeed the enveloping algebra of its associated Lie algebra
${\mathcal L}^1$ is
the dual of the Connes-Kreimer Hopf algebra describing renormalisation
theory in \cite{CK98}. Foissy proved that the Lie algebra ${\mathcal
  L}^1$ is a free Lie algebra in \cite{F02}. This results generalize
easily to the following statement: the Lie algebra associated to a
free pre-Lie algebra spanned by a vector space $V$ is a free Lie
algebra spanned by a vector space $W$. However, the proof of Foissy
doesn't give an explicit description of $W$. 

In the language of species \cite{Aguiar,Joyal86,L06}, a symmetric operad is a monoid in the category of species with respect to the composition product.
Foissy's result suggests that at the
level of species, there exists a specie $\mathcal F$ such that
$$\calT=\Lie\circ\calF.$$
Indeed, this is proved by Chapoton in \cite{Chapo07} where he describes $\calF$
by the way of symmetric functions. 

In this paper we aim to give 
an explicit basis for  $\calF$. This will allows us in Section~\ref{S:LieF} and Section~\ref{S:LieFsp} to give a new proof of Chapoton's result. We also show in Section~\ref{S:F} that $\calF$ is a sub nonsymmetric operad of the pre-Lie operad $\calT$.
Before all this, we recall some basic fact in Section \ref{S:LieT}.
 
  
\section{The pre-Lie operad and rooted trees}\label{S:LieT}

 We first recall the definition of the pre-Lie operad based
  on labelled rooted trees as in \cite{CL01}. For $n\in \mathbb N^*$,
the set $\{1,\ldots,n\}$ is denoted by $[n]$ and $[0]$ denotes the empty set.
The symmetric group on $k$ letters is denoted by $S_k$.

Recall that a specie is a contravariant functor from the category of
finite sets $\hbox{\bf
  Set}^\times$ and bijections to the category of
finite dimensional vector spaces $\hbox{\bf Vect}$.
Following Joyal in \cite{Joyal86}, a specie is equivalent to an
$\mathbb S$-module, that is, a collection of vector spaces
$(V_n)_{n\geq 0}$ with a right action of $S_n$ on $V_n$.

 Given two species $\calA,\calB\colon\hbox{\bf Set}^\times\to\hbox{\bf Vect}$ we have the product
    \begin{equation}\label{eq:prod}
    \calA\bullet \calB [S] = \bigoplus_{I+J=S} \calA[I]\otimes \calB[J]\,,
    \end{equation}
where $I+J$ denotes the disjoint union of the sets $I$ and $J$. 
We have the composition of species defined by
 \begin{equation*}
  \calA\circ\calB [S]= \bigoplus_{k\geq 0} \ \calA[k] \otimes_{S_k}
  (\calB^{\bullet k} [S])\,.
  \end{equation*}

If $\calB[\emptyset]=0$ the composition of species has the form
  \begin{equation}\label{eq:comp}
  \calA\circ\calB [S]= \bigoplus_{\Phi\parti S} \ \calA[\Phi] \otimes \left({\bigotimes_{\phi\in\Phi}} \calB[\phi]\right)\,,
  \end{equation}
where $\Phi\parti S$ denotes that $\Phi$ is a set partition of $S$.  

A symmetric operad is a monoid in the category of species with respect
to the composition product. A (twisted) algebra $A$ over an operad
$\mathcal P$ is a specie together with an evaluation product
$$\mu_A:\mathcal P\circ A \rightarrow A$$
satisfying the usual condition (see \cite{Aguiar} or \cite{L06} for more details). Note that 
what is called an algebra over an operad is usually a vector
space considered as a specie which is always zero except on the emptyset.
The terminology ``twisted'' emphasizes the fact that we generalize the usual
definition to any specie. For instance $\mathcal P$ is the free
twisted $\mathcal P$-algebra generated by the unit $I$ for the
composition product, whereas $\oplus_{n\geq 0} \mathcal P[n]/S_n$ is
the free $\mathcal P$-algebra generated by a 1-dimensional vector
space.

\smallskip

  Given a finite set $S$ of cardinality $n$ let $\calT[S]$ be the
  vector space freely generated by the labelled rooted  trees on $n$
  vertices with distinct label chosen in $S$. For $n=0$ we set
$T[\emptyset]=0$. This gives us a specie.

 \begin{exa}\label{e:T3} The space $\calT[\{1,2,3\}]$ is the linear span of the following trees:
   $$
 \raise -10pt
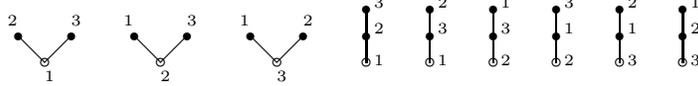
\hbox{ \begin{picture}(40,30)
       \put(10,6){\circle{3}}   \put(0,16){\circle*{3}}   \put(20,16){\circle*{3}}
       \put(10,6){\line(-1,1){10}}    \put(10,6){\line(1,1){10}}
       \put(10,-1){$\scriptscriptstyle 1$}   \put(-4,20){$\scriptscriptstyle 2$}   \put(20,20){$\scriptscriptstyle 3$}  
      \end{picture}} 
 \raise -10pt\hbox{  \begin{picture}(40,30)
       \put(10,6){\circle{3}} \put(0,16){\circle*{3}} \put(20,16){\circle*{3}}
       \put(10,6){\line(-1,1){10}}  \put(10,6){\line(1,1){10}}
       \put(10,-1){$\scriptscriptstyle 2$}   \put(-4,20){$\scriptscriptstyle 1$} \put(20,20){$\scriptscriptstyle 3$}  
      \end{picture}} 
 \raise -10pt\hbox{ \begin{picture}(40,30)
       \put(10,6){\circle{3}} \put(0,16){\circle*{3}}  \put(20,16){\circle*{3}}
       \put(10,6){\line(-1,1){10}}  \put(10,6){\line(1,1){10}}
       \put(10,-1){$\scriptscriptstyle 3$}   \put(-4,20){$\scriptscriptstyle 1$} \put(20,20){$\scriptscriptstyle 2$}  
      \end{picture}} 
 \raise -10pt\hbox{ \begin{picture}(20,30)
       \put(0,6){\circle{3}} \put(0,16){\circle*{3}} \put(0,26){\circle*{3}}
       \put(0,6){\line(0,1){20}}
       \put(3,5){$\scriptscriptstyle 1$} \put(3,17){$\scriptscriptstyle 2$} \put(3,27){$\scriptscriptstyle 3$}  
      \end{picture}} 
 \raise -10pt\hbox{ \begin{picture}(20,30)
       \put(0,6){\circle{3}}  \put(0,16){\circle*{3}} \put(0,26){\circle*{3}}
       \put(0,6){\line(0,1){20}}
       \put(3,5){$\scriptscriptstyle 1$}  \put(3,17){$\scriptscriptstyle 3$}  \put(3,27){$\scriptscriptstyle 2$}  
      \end{picture}} 
 \raise -10pt\hbox{ \begin{picture}(20,30)
       \put(0,6){\circle{3}}  \put(0,16){\circle*{3}} \put(0,26){\circle*{3}}
       \put(0,6){\line(0,1){20}} 
        \put(3,5){$\scriptscriptstyle 2$}  \put(3,17){$\scriptscriptstyle 3$} \put(3,27){$\scriptscriptstyle 1$}  
      \end{picture}} 
 \raise -10pt\hbox{ \begin{picture}(20,30)
       \put(0,6){\circle{3}}  \put(0,16){\circle*{3}} \put(0,26){\circle*{3}}
       \put(0,6){\line(0,1){20}}
       \put(3,5){$\scriptscriptstyle 2$}  \put(3,17){$\scriptscriptstyle 1$}  \put(3,27){$\scriptscriptstyle 3$}  
      \end{picture}} 
 \raise -10pt\hbox{ \begin{picture}(20,30)
       \put(0,6){\circle{3}} \put(0,16){\circle*{3}}  \put(0,26){\circle*{3}}
       \put(0,6){\line(0,1){20}}
       \put(3,5){$\scriptscriptstyle 3$}  \put(3,17){$\scriptscriptstyle 1$}  \put(3,27){$\scriptscriptstyle 2$}  
      \end{picture}} 
 \raise -10pt\hbox{ \begin{picture}(20,30)
       \put(0,6){\circle{3}}  \put(0,16){\circle*{3}} \put(0,26){\circle*{3}}
       \put(0,6){\line(0,1){20}}
       \put(3,5){$\scriptscriptstyle 3$}   \put(3,17){$\scriptscriptstyle 2$}   \put(3,27){$\scriptscriptstyle 1$}  
      \end{picture}} 
 $$
In general there are $n^{n-1}$ such trees on a set of cardinality $n$ (see \cite{BLL} for more details). 
  \end{exa}
  
\begin{theo}{\cite[theorem 1.9]{CL01}} The specie $\mathcal T$ forms
  an operad. Algebras over this operad are pre-Lie algebras, that is,
  vector spaces $L$ together with a product $*$ satisfying the
  relation
$$(x* y)* z-x*(y* z)=(x* z)* y-x*(z* y), \ \forall x,y,z\in L.$$
\end{theo}
As a consequence $\mathcal T$ is the free twisted pre-Lie algebra
generated by $I$. The twisted pre-Lie product is described as follows.

\begin{defi}\label{D:pre-Lieproduct}
Given two disjoint sets $I,J$ and two trees $T\in\calT[I]$ and $Y\in\calT[J]$ we define
$$T*Y=
\def\objectstyle{\scriptstyle}
\def\labelstyle{\scriptstyle}
\sum_{t\in Vert(T)}
\vcenter{\xymatrix@-1.5pc{
*++[o][F-]{Y}\ar@{-*{\bullet}}[d]_>>{t}\\
*++[o][F-]{T}}}$$
where the sum is over all possible ways of {\sl grafting} the root of
the tree $Y$ on a vertex $t$ of $T$. The root of $T*Y$ is the one of $T$.
\end{defi} 

Since any pre-Lie algebra $L$ gives rise to a Lie algebra whose bracket
is defined by
$[x,y]=x*y-y*x$ there is a morphism of operads 
$$\Lie \rightarrow \mathcal T.$$
Note that this morphism is injective: an associative algebra is
obviously a pre-Lie algebra and the composition of morphisms of operads
$$\Lie\rightarrow \mathcal T \rightarrow {\mathcal A}s$$
is the usual injective morphism from $\Lie$ to ${\mathcal A}s$.
As a consequence the specie $\mathcal T$ is a twisted Lie algebra, that
is a Lie monoid in the category of specie. It is endowed with the
following Lie
bracket $[\,\,,\,]\colon\calT\bullet\calT\to\calT$: given two disjoint 
sets $I,J$ and two trees $T\in\calT[I]$ and $Y\in\calT[J]$ we define

\begin{equation}\label{D:bracket}
[T,Y]=T*Y-Y*T=
\def\objectstyle{\scriptstyle}
\def\labelstyle{\scriptstyle}
\sum_{t\in Vert(T)}
\vcenter{\xymatrix@-1.5pc{
*++[o][F-]{Y}\ar@{-*{\bullet}}[d]_>>{t}\\
*++[o][F-]{T}}}
-\sum_{s\in Vert(Y)}
\vcenter{\xymatrix@-1.5pc{
*++[o][F-]{T}\ar@{-*{\bullet}}[d]_>>{s}\\
*++[o][F-]{Y}}}
\end{equation}

 \begin{exa} For 
 $T= \raise -10pt\hbox{ \begin{picture}(25,25)
       \put(10,6){\circle{3}}   \put(0,16){\circle*{3}}   \put(20,16){\circle*{3}}
       \put(10,6){\line(-1,1){10}}    \put(10,6){\line(1,1){10}}
       \put(10,-1){$\scriptscriptstyle 3$}   \put(-4,20){$\scriptscriptstyle 1$}   \put(20,20){$\scriptscriptstyle 4$}  
       \end{picture}} 
       \in\calT[\{1,3,4\}]$ 
 and 
 $Y=\raise -5pt\hbox{ \begin{picture}(10,5)
        \put(0,6){\circle{3}}     \put(3,5){$\scriptscriptstyle 2$}    \end{picture}}
        \in\calT[\{2\}]$ 
we have that
   $$
   [\raise -10pt\hbox{ \begin{picture}(25,25)
       \put(10,6){\circle{3}}   \put(0,16){\circle*{3}}   \put(20,16){\circle*{3}}
       \put(10,6){\line(-1,1){10}}    \put(10,6){\line(1,1){10}}
       \put(10,-1){$\scriptscriptstyle 3$}   \put(-4,20){$\scriptscriptstyle 1$}  
        \put(20,20){$\scriptscriptstyle 4$}  
       \end{picture}} \,,
    \raise -5pt\hbox{ \begin{picture}(10,5)
        \put(0,6){\circle{3}}     \put(3,5){$\scriptscriptstyle 2$}    \end{picture}}]
    \quad = \quad   
    \raise -10pt\hbox{ \begin{picture}(25,30)
       \put(10,6){\circle{3}}   \put(0,16){\circle*{3}}   \put(20,16){\circle*{3}}      \put(0,26){\circle*{3}}
       \put(10,6){\line(-1,1){10}}    \put(10,6){\line(1,1){10}}                           \put(0,16){\line(0,1){10}}        
       \put(10,-1){$\scriptscriptstyle 3$}   \put(-5,18){$\scriptscriptstyle 1$}  
        \put(20,20){$\scriptscriptstyle 4$}                                           \put(-4,30){$\scriptscriptstyle 2$} 
       \end{picture}}
    \quad + \quad   
    \raise -10pt\hbox{ \begin{picture}(25,30)
       \put(10,6){\circle{3}}   \put(0,16){\circle*{3}}   \put(20,16){\circle*{3}}      \put(10,16){\circle*{3}}
       \put(10,6){\line(-1,1){10}}    \put(10,6){\line(1,1){10}}                             \put(10,6){\line(0,1){10}}        
       \put(10,-1){$\scriptscriptstyle 3$}   \put(-4,20){$\scriptscriptstyle 1$}  
        \put(20,20){$\scriptscriptstyle 4$}                                           \put(7,20){$\scriptscriptstyle 2$} 
       \end{picture}}
     \quad + \quad   
     \raise -10pt\hbox{ \begin{picture}(25,30)
       \put(10,6){\circle{3}}   \put(0,16){\circle*{3}}   \put(20,16){\circle*{3}}      \put(20,26){\circle*{3}}
       \put(10,6){\line(-1,1){10}}    \put(10,6){\line(1,1){10}}                           \put(20,16){\line(0,1){10}}        
       \put(10,-1){$\scriptscriptstyle 3$}   \put(-5,20){$\scriptscriptstyle 1$}  
        \put(22,18){$\scriptscriptstyle 4$}                                           \put(20,30){$\scriptscriptstyle 2$} 
       \end{picture}}
     \quad - \quad   
     \raise -10pt\hbox{ \begin{picture}(25,30)
       \put(10,16){\circle*{3}}   \put(0,26){\circle*{3}}   \put(20,26){\circle*{3}}      \put(10,6){\circle{3}}
       \put(10,16){\line(-1,1){10}}    \put(10,16){\line(1,1){10}}                           \put(10,6){\line(0,1){10}}        
       \put(12,13){$\scriptscriptstyle 3$}   \put(-5,30){$\scriptscriptstyle 1$}  
        \put(22,28){$\scriptscriptstyle 4$}                                           \put(10,-1){$\scriptscriptstyle 2$} 
       \end{picture}}\,.
 $$
\end{exa}

\medskip

As we mentioned in the Introduction, we shall now describe explicitly  a specie $\mathcal F$ such that
$\calT=\Lie\circ\calF.$

\section{$\calT[S] = \Lie \circ \calF[S]$ as vector spaces}\label{S:LieF}

In this section, we show an auxiliary result relating  $\calT$ to a free Lie algebra over rooted trees that are increasing in the first level. We give an explicit  isomorphism using bases. This has the advantage to be explicit but it is not natural. In the next section we will induce an action of the symmetric groups on both side hence giving an identity of species.

Given a finite set $S$ and a linear order on $S$, let $\calF[S]$ be the vector space spanned by the basis of $S$-labelled rooted trees that are increasing at the first level. That is the trees such that the labels increase from the root to the adjacent vertices and no other condition on the other labels. Also, we let $\calF[\emptyset]=0$. At this point, 
$\calF$ is not a specie as it depends on an order on $S$. 
 We will turn this into a specie in the next section.

 \begin{exa}\label{e:F3} The space $\calF[\{1,2,3\}]$ with the natural order on  $\{1,2,3\}$ has basis given by the following trees:
   $$
 \raise -10pt
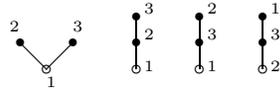
\hbox{ \begin{picture}(40,30)
       \put(10,6){\circle{3}}   \put(0,16){\circle*{3}}   \put(20,16){\circle*{3}}
       \put(10,6){\line(-1,1){10}}    \put(10,6){\line(1,1){10}}
       \put(10,-1){$\scriptscriptstyle 1$}   \put(-4,20){$\scriptscriptstyle 2$}   \put(20,20){$\scriptscriptstyle 3$}  
      \end{picture}} 
 \raise -10pt\hbox{ \begin{picture}(20,30)
       \put(0,6){\circle{3}} \put(0,16){\circle*{3}} \put(0,26){\circle*{3}}
       \put(0,6){\line(0,1){20}}
       \put(3,5){$\scriptscriptstyle 1$} \put(3,17){$\scriptscriptstyle 2$} \put(3,27){$\scriptscriptstyle 3$}  
      \end{picture}} 
 \raise -10pt\hbox{ \begin{picture}(20,30)
       \put(0,6){\circle{3}}  \put(0,16){\circle*{3}} \put(0,26){\circle*{3}}
       \put(0,6){\line(0,1){20}}
       \put(3,5){$\scriptscriptstyle 1$}  \put(3,17){$\scriptscriptstyle 3$}  \put(3,27){$\scriptscriptstyle 2$}  
      \end{picture}} 
 \raise -10pt\hbox{ \begin{picture}(20,30)
       \put(0,6){\circle{3}}  \put(0,16){\circle*{3}} \put(0,26){\circle*{3}}
       \put(0,6){\line(0,1){20}} 
        \put(3,5){$\scriptscriptstyle 2$}  \put(3,17){$\scriptscriptstyle 3$} \put(3,27){$\scriptscriptstyle 1$}  
      \end{picture}} 
 $$ 
 In general there are $(n-1)^{n-1}$ such trees (see e.g. \cite{CDG00} for more details)
\end{exa}

For our next result, we also need to consider $\Lie[S]$ as the vector space of multilinear brackets of  degree $|S|$. That is the vector space spanned by all brackets of the elements of $S$ (without repetition) modulo the antisymmetry relation and Jacobi identity. It is easy to check that this construction is functorial. We then have that $\Lie$ is a specie. In fact, 
$\Lie$ is an operad and algebras over this operad are the classical Lie algebras. (see e.g. \cite{GK94}).

 \begin{exa}\label{ex:lie} The space $\Lie[\{1,2,3\}]$ is the linear span of the following brackets:
   $$[[1,2],3],\  [[1,3],2],\  [[2,1],3],\  [[2,3],1],\  [[3,1],2],\  [[3,2],1],$$
   $$[1,[2,3]],\  [1,[3,2]],\  [2,[1,3]],\  [2,[3,1]],\  [3,[1,2]],\  [3,[2,1]].$$ 
As we will see below, it is well known that this space has dimension
equal to two and that a basis is given 
by $\{[[3,1],2], [3,[2,1]]\}$.  In general there are $(n-1)!$  linearly independent brackets \cite{Reut}.
\end{exa}

If we are given a linear order on $S$ we can construct an explicit basis of $\Lie[S]$. This is the classical Lyndon basis of $\Lie$ (see \cite{Reut}).
More precisely, $\Lie[S]$ has basis given by the Lyndon permutations with Lyndon bracketing. For our purpose we use the reverse lexicographic order to produce the following basis of $\Lie[S]$. 
Let $S=\{a<b<\cdots<y<z\}$. 
A {\sl Lyndon} permutation  $\sigma\colon S\to S$  is a permutation such that $\sigma(a)=z$. The Lyndon bracketing $\sb[\sigma]$ of $\sigma$ is defined recursively. We write $\sigma=\big(\sigma(a),\sigma(b),\ldots,\sigma(z)\big)$ as the list of its values. 
If $S=\{a\}$ is of cardinality 1, then define $\sb[\sigma(a)]=a$. If $|S|>1$, let $k\in S$ be such that $\sigma(k)=y$ the second largest value of $S$, then define
 $$\sb[\sigma(a),\ldots,\sigma(j),\sigma(k),\ldots,\sigma(z)] =
     \big[\sb[\sigma(a),\ldots,\sigma(j)],\ \sb[\sigma(k),\ldots,\sigma(z)]\big].$$
     A basis of $\Lie[S]$ is given by the set $\{\sb[\sigma]:\sigma \hbox{ is Lyndon}\}$
In the Example~\ref{ex:lie} we have that $(3,1,2)$ and $(3,2,1)$ are the only two Lyndon permutations and $\sb[3,1,2]=[\sb[3,1],\sb[2]]=[[\sb[3],\sb[1]],2]=[[3,1],2]$. Similarly $\sb[3,2,1]=[3,[2,1]]$. 

Even though $\calF$ is not a specie we can still define $\Lie\circ \calF$. Let $S$ be a finite set with a linear order. We define
$$  \Lie\circ\calF [S]= \bigoplus_{\Phi\parti S} \ \Lie[\Phi] \otimes \left({\bigotimes_{\phi\in\Phi}} \calF[\phi]\right)\,,$$
where for $\Phi\parti S$ we induce a linear order on each part $\phi\in\Phi$ from the linear order on $S$.  

\begin{theo} \label{th:TLieF1}
    $\calT[S] = \Lie\circ \calF[S]$ as vector spaces. 
\end{theo}

\begin{proof}
Given a linear order on a finite set $S$, we construct a linear isomorphism between $\calT[S]$ and $\Lie\circ \calF[S]$. By definition, $\Lie\circ \calF[S]$ is any bracketing of trees of type $\calF$ such that the disjoint union of all the labels is $S$. Since $\calT$ is a Lie monoid there is a natural map $\Xi\colon\Lie\circ \calF[S]\to\calT[S]$. We need to show that this map is injective and surjective.

Assume that we have a finite set $S$ and a linear order on $S$. For $\Phi=\{\phi_1,\phi_2,\ldots,\phi_\ell\}\parti S$ we have that each part $\phi_i$ is also ordered. We can then order any set of trees $\{T_i : T_i\in\calF[\phi_i],\ 1\le i\le \ell\}$ using the roots of the trees. It follows  that a basis for  $\Lie\circ \calF[S]$ is given by 
  $$\big\{\sb[T_{\sigma(1)}\cdots T_{\sigma(\ell)}]:\substack{\Phi=\{\phi_1,\ldots,\phi_\ell\}\parti S,\
               \sigma\colon[\ell]\to[\ell],\ T_i\in\calF[\phi_i],\ T_{\sigma(1)} \hbox{ \small has the largest root}}
      \big\}\,.$$
 To complete the proof, we need to show that
\begin{equation}\label{eq:basis}
\big\{\Xi(\sb[T_{\sigma(1)}\cdots T_{\sigma(\ell)}]):\substack{\Phi=\{\phi_1,\ldots,\phi_\ell\}\parti S,\
               \sigma\colon[\ell]\to[\ell],\ T_i\in\calF[\phi_i],\ T_{\sigma(1)} \hbox{ \small has the largest root}}
      \big\}\,.
      \end{equation}
is a basis of $\calT[S]$. Using the order on $S$, we introduce a grading on the basis of labelled rooted  trees of $\calT[S]$ and show that there exists a triangularity relation between the basis in (\ref{eq:basis}) and the basis of labelled rooted trees. We say that a tree $T\in\calT[S]$ is of degree $d$ if the maximal decreasing connected subtree of $T$ from the root has $d$ vertices. For any tree $T\in\calT[S]$ we denote by $MD(T)$ its maximal decreasing connected subtree from the root.
For example consider
   $$
T_1=\quad \raise -10pt\hbox{ \begin{picture}(40,30)
       \put(10,6){\circle{3}}   \put(0,16){\circle*{3}}   \put(20,16){\circle*{3}}  \put(10,16){\circle*{3}}
        \put(20,26){\circle*{3}} \put(0,26){\circle*{3}} \put(-10,26){\circle*{3}}
       \put(10,6){\line(-1,1){20}}    \put(10,6){\line(1,1){10}} 
         \put(10,6){\line(0,1){10}}    \put(0,16){\line(0,1){10}}    \put(20,16){\line(0,1){10}} 
        \put(10,-1){$\scriptscriptstyle 5$}   \put(-4,10){$\scriptscriptstyle 2$}   \put(22,18){$\scriptscriptstyle 6$}  
       \put(-14,30){$\scriptscriptstyle 7$}   \put(-2,30){$\scriptscriptstyle 1$}  
        \put(8,20){$\scriptscriptstyle 3$}   \put(20,30){$\scriptscriptstyle 4$}  
      \end{picture}} 
      \qquad\hbox{and}\qquad
T_2=      \raise -10pt\hbox{ \begin{picture}(25,30)
       \put(10,16){\circle*{3}}   \put(0,26){\circle*{3}}   \put(20,26){\circle*{3}}      \put(10,6){\circle{3}}
       \put(10,16){\line(-1,1){10}}    \put(10,16){\line(1,1){10}}                           \put(10,6){\line(0,1){10}}        
       \put(12,13){$\scriptscriptstyle 3$}   \put(-5,30){$\scriptscriptstyle 1$}  
        \put(22,28){$\scriptscriptstyle 4$}                                           \put(10,-1){$\scriptscriptstyle 2$} 
       \end{picture}}\,.
 $$ 
$MD(T_1)$ is build with the vertices labelled $\{5,3,2,1\}$ and for $MD(T_2)$ we use only $\{2\}$. Hence $T_1$ is of degree 4 and $T_2$ is of degree 1. Remark that $T\in \calF[S]$ if and only if the degree of $T$ is 1. 

Given a set partition $\Phi=\{\phi_1,\ldots,\phi_\ell\}\parti S$, a permutation $\sigma\colon[\ell]\to[\ell]$, a family of trees  $\{T_i:T_i\in\calF[\phi_i],\ T_{\sigma(1)} \hbox{ \small has the largest root}\}$, we claim that in the expansion of $\Xi(\sb[T_{\sigma(1)}\cdots T_{\sigma(\ell)}])$ there is a unique tree  of maximal degree $\ell$ (with coefficient 1). Furthermore, the correspondence from $\Xi(\sb[T_{\sigma(1)}\cdots T_{\sigma(\ell)}])$ to its maximal degree term $T$ is such that $MD(T)$ is formed from the vertices labelled by the labels of the roots of $T_1, T_2, \ldots, T_\ell$. In fact  the maximal decreasing subtree  of any tree in the expansion of $\Xi(\sb[T_{\sigma(1)}\cdots T_{\sigma(\ell)}])$ is formed from the vertices labelled by a subset of the labels of the roots of $T_1, T_2, \ldots, T_\ell$.

We proceed by induction on $\ell$. For $\ell=1$ we have $\Xi(\sb[T_1])=T_1$ a unique tree of degree 1. For $\ell=2$ we are given two trees of degree 1:
$$
T_1=\quad \raise -10pt\hbox{ \begin{picture}(40,30)
       \put(10,6){\circle{3}}   \put(0,16){\circle*{3}}   \put(-10,16){\circle*{3}}  \put(30,16){\circle*{3}}
       \put(10,6){\line(-2,1){20}}    \put(10,6){\line(-1,1){10}}   \put(10,6){\line(2,1){20}}     
       \put(10,-1){$\scriptscriptstyle b$}   \put(10,20){$\scriptscriptstyle \cdots$}   
       \put(-14,20){$\scriptscriptstyle Y_1$}   \put(-2,20){$\scriptscriptstyle Y_2$}  
        \put(28,20){$\scriptscriptstyle Y_r$}    
      \end{picture}} 
      \qquad\hbox{and}\qquad
T_2= \quad \raise -10pt\hbox{ \begin{picture}(40,30)
       \put(10,6){\circle{3}}   \put(0,16){\circle*{3}}   \put(-10,16){\circle*{3}}  \put(30,16){\circle*{3}}
       \put(10,6){\line(-2,1){20}}    \put(10,6){\line(-1,1){10}}   \put(10,6){\line(2,1){20}}     
       \put(10,-1){$\scriptscriptstyle a$}   \put(10,20){$\scriptscriptstyle \cdots$}   
       \put(-14,20){$\scriptscriptstyle X_1$}   \put(-2,20){$\scriptscriptstyle X_2$}  
        \put(28,20){$\scriptscriptstyle X_k$}    
      \end{picture}} \,,
$$
 where $T_i\in\calF[\phi_i]$. This implies that the roots of each $Y_j$ is strictly greater than $b$ and the roots of each $X_j$ is strictly greater than $a$. We assume without lost of generality that $b>a$. When we expand $\Xi(\sb[T_1T_2])=[T_1,T_2]$ we obtain
$$
 \raise -10pt\hbox{ \begin{picture}(50,50)
       \put(10,6){\circle{3}}   \put(0,16){\circle*{3}}   \put(-10,16){\circle*{3}}  \put(20,16){\circle*{3}}
       \put(10,6){\line(-2,1){20}}    \put(10,6){\line(-1,1){10}}   \put(10,6){\line(1,1){10}}     
       \put(10,-1){$\scriptscriptstyle b$}   \put(10,20){$\scriptscriptstyle \cdots$}   
       \put(-14,20){$\scriptscriptstyle Y_1$}   \put(-2,20){$\scriptscriptstyle Y_2$}  
        \put(20,20){$\scriptscriptstyle Y_r$}    
    \put(40,16){\circle*{3}}   \put(40,26){\circle*{3}}   \put(30,26){\circle*{3}}  \put(60,26){\circle*{3}}
       \put(40,16){\line(-1,1){10}}    \put(40,16){\line(0,1){10}}   \put(40,16){\line(2,1){20}}     
       \put(40,9){$\scriptscriptstyle a$}   \put(50,30){$\scriptscriptstyle \cdots$}   
       \put(26,30){$\scriptscriptstyle X_1$}   \put(38,30){$\scriptscriptstyle X_2$}  
        \put(58,30){$\scriptscriptstyle X_k$}     \put(10,6){\line(3,1){30}} 
      \end{picture}} 
      \quad-\quad
 \raise -10pt\hbox{ \begin{picture}(50,50)
       \put(10,6){\circle{3}}   \put(0,16){\circle*{3}}   \put(-10,16){\circle*{3}}  \put(20,16){\circle*{3}}
       \put(10,6){\line(-2,1){20}}    \put(10,6){\line(-1,1){10}}   \put(10,6){\line(1,1){10}}     
       \put(10,-1){$\scriptscriptstyle a$}   \put(10,20){$\scriptscriptstyle \cdots$}   
       \put(-14,20){$\scriptscriptstyle X_1$}   \put(-2,20){$\scriptscriptstyle X_2$}  
        \put(20,20){$\scriptscriptstyle X_k$}    
    \put(40,16){\circle*{3}}   \put(40,26){\circle*{3}}   \put(30,26){\circle*{3}}  \put(60,26){\circle*{3}}
       \put(40,16){\line(-1,1){10}}    \put(40,16){\line(0,1){10}}   \put(40,16){\line(2,1){20}}     
       \put(40,9){$\scriptscriptstyle b$}   \put(50,30){$\scriptscriptstyle \cdots$}   
       \put(26,30){$\scriptscriptstyle Y_1$}   \put(38,30){$\scriptscriptstyle Y_2$}  
        \put(58,30){$\scriptscriptstyle Y_r$}     \put(10,6){\line(3,1){30}} 
      \end{picture}} 
  \quad+\ \sum\quad
       \raise -10pt\hbox{ \begin{picture}(40,60)
       \put(10,6){\circle{3}}   \put(0,16){\circle*{3}}   \put(-10,16){\circle*{3}}  \put(30,16){\circle*{3}}
       \put(10,6){\line(-2,1){20}}    \put(10,6){\line(-1,1){10}}   \put(10,6){\line(2,1){20}}     
       \put(10,-1){$\scriptscriptstyle b$}   \put(10,20){$\scriptscriptstyle \cdots$}   
       \put(-14,20){$\scriptscriptstyle Y_1$}   \put(-2,20){$\scriptscriptstyle Y_2$}  
        \put(28,20){$\scriptscriptstyle Y_r$}    
    \put(20,36){\circle*{3}}   \put(10,46){\circle*{3}}   \put(0,46){\circle*{3}}  \put(40,46){\circle*{3}}
       \put(20,36){\line(-2,1){20}}    \put(20,36){\line(-1,1){10}}   \put(20,36){\line(2,1){20}}     
       \put(20,29){$\scriptscriptstyle a$}   \put(20,50){$\scriptscriptstyle \cdots$}   
       \put(-4,50){$\scriptscriptstyle X_1$}   \put(8,50){$\scriptscriptstyle X_2$}  
        \put(38,50){$\scriptscriptstyle X_k$}     \put(20,36){\line(-1,-2){8}}  
      \end{picture}} 
  \quad-\ \sum\quad
\raise -10pt\hbox{ \begin{picture}(40,60)
       \put(20,36){\circle{3}}   \put(10,46){\circle*{3}}   \put(0,46){\circle*{3}}  \put(40,46){\circle*{3}}
       \put(20,36){\line(-2,1){20}}    \put(20,36){\line(-1,1){10}}   \put(20,36){\line(2,1){20}}     
       \put(20,29){$\scriptscriptstyle b$}   \put(20,50){$\scriptscriptstyle \cdots$}   
       \put(-4,50){$\scriptscriptstyle Y_1$}   \put(8,50){$\scriptscriptstyle Y_2$}  
        \put(38,50){$\scriptscriptstyle Y_r$}    
     \put(10,6){\circle{3}}   \put(0,16){\circle*{3}}   \put(-10,16){\circle*{3}}  \put(30,16){\circle*{3}}
       \put(10,6){\line(-2,1){20}}    \put(10,6){\line(-1,1){10}}   \put(10,6){\line(2,1){20}}     
       \put(10,-1){$\scriptscriptstyle a$}   \put(10,20){$\scriptscriptstyle \cdots$}   
       \put(-14,20){$\scriptscriptstyle X_1$}   \put(-2,20){$\scriptscriptstyle X_2$}  
        \put(28,20){$\scriptscriptstyle X_k$}    \put(20,36){\line(-1,-2){8}}  
      \end{picture}} \,.
$$
The first term is of degree 2 and its maximal decreasing subtree is build from $\{b,a\}$ the roots of $T_1$ and $T_2$. All the other trees in this expansion are of degree 1 and their maximal decreasing subtrees are labelled either by $a$ or by $b$ . 

We now assume that $\ell>2$. To compute $\Xi(\sb[T_{\sigma(1)}\cdots T_{\sigma(\ell)}])$, let  $b_1, b_2, \ldots, b_\ell$ be the roots of $T_{\sigma(1)}, T_{\sigma(2)}, \ldots T_{\sigma(\ell)}$ respectively. By construction we have that $b_1=\max(b_1, b_2, \ldots, b_\ell)$. Let $b_k=\max(b_2, \ldots, b_\ell)$. That is $b_k$ is the second largest root and $k>1$. The Lyndon factorization gives us that 
 $$\Xi(\sb[T_{\sigma(1)}\cdots T_{\sigma(\ell)}])=
     \big[ \Xi(\sb[T_{\sigma(1)}\cdots T_{\sigma(k-1)}]),\Xi(\sb[T_{\sigma(k)}\cdots T_{\sigma(\ell)}])\big].$$
By induction hypothesis we have that 
     $$\Xi(\sb[T_{\sigma(1)}\cdots T_{\sigma(k-1)}])= Y_0 +\sum_i c_i Y_i$$
where $Y_0$ is of degree $k-1$ and $MD(Y_0)$ is formed with vertices labelled by $\{b_1,\ldots,b_{k-1}\}$
and the trees $Y_i$ ($i\ne 0$) are of degree $<k-1$ where $MD(Y_i)$ are formed with vertices labelled by a subset of $\{b_1,\ldots,b_{k-1}\}$. Similarly, 
     $$\Xi(\sb[T_{\sigma(k)}\cdots T_{\sigma(\ell)}])=X_0 + \sum_j d_j X_j$$
where $X_0$ is of degree $\ell-k+1$ and $MD(X_0)$ is formed with vertices labelled by 
$\{b_k,\ldots,b_{\ell}\}$
and the trees $X_j$ ($j\ne 0$) are of degree $<\ell-k+1$ where $MD(X_j)$ are formed with vertices labelled by a subset of $\{b_k,\ldots,b_{\ell}\}$ .  
The largest degree term in $[Y_i,X_j]$ must be obtained by either
grafting $MD(Y_i)$ in $MD(X_j)$, or by grafting$MD(X_j)$ in
$MD(Y_i)$. Hence the 
largest degree term in $[Y_i,X_j]$ is of degree at most $\deg(Y_i)+\deg(X_j)$. Hence it is sufficient to concentrate our attention on $[Y_0,X_0]$. In this case, recall that $b_1$ is the largest value, so it must be the root of $MD(Y_0)$. Similarly $b_k$ is the root of $MD(X_0)$. We can get a tree of degree $\ell$ by grafting $X_0$ at the root of $Y_0$. If we graft $X_0$ anywhere else in $Y_0$ we get a tree of degree strictly smaller.
In fact, since $b_k>\max(b_2,...,b_{k-1})$, if we graft $X_0$ on $MD(Y_0)$ (not at the root) or anywhere else, we get a tree of degree equal to $\deg(Y_0)=k-1<\ell$. On the other hand, since $b_1$ is maximal, if we graft $Y_0$ in $X_0$ we always get a tree of degree equal to $\deg(X_0)=\ell-k+1<\ell$.

We now remark that  $MD(Z)$ of any term $Z$ in the expansion of $[Y_i,X_j]$, is either $MD(Y_i)$,$MD(X_j)$, the grafting of $MD(Y_i)$ in $MD(X_j)$ or the grafting of $MD(X_j)$ in $MD(Y_i)$. In all cases, the vertices of $MD(Z)$ are labelled by a subset of $\{b_1,\ldots,b_\ell\}$ and this conclude the induction.

To conclude the triangularity relation we need to show that for any tree $T\in \calT[S]$ there is a basis element in the basis (\ref{eq:basis}) with $T$ as its leading degree term. For this we proceed by induction on the degree of $T$. Our hypothesis is that for any tree $T\in \calT[S]$ we can find
 a set partition $\Phi=\{\phi_1,\ldots,\phi_\ell\}\parti S$, a permutation $\sigma\colon[\ell]\to[\ell]$ and a family of trees  $\{T_i:T_i\in\calF[\phi_i],\ T_{\sigma(1)} \hbox{ \small has the largest root}\}$, such that $T$ is the leading term of  $\Xi(\sb[T_{\sigma(1)}\cdots T_{\sigma(\ell)}])$. Furthermore, $MD(T)$ is the subtree formed with the vertices labelled by  labels of the roots of $T_1,\ldots, T_\ell$.

 If $T$ is of degree 1, then $T=\Xi(\sb[T])$ and $MD(T)$ is a single vertex. If $deg(T)>1$, then $T$ is of the form
\begin{equation}\label{eq:T}
T\quad=\quad
 \raise -10pt\hbox{ \begin{picture}(50,50)
       \put(10,6){\circle{3}}   \put(0,16){\circle*{3}}   \put(-10,16){\circle*{3}}  \put(20,16){\circle*{3}}
       \put(10,6){\line(-2,1){20}}    \put(10,6){\line(-1,1){10}}   \put(10,6){\line(1,1){10}}     
       \put(10,-1){$\scriptscriptstyle b$}   \put(10,20){$\scriptscriptstyle \cdots$}   
       \put(-14,20){$\scriptscriptstyle Y_1$}   \put(-2,20){$\scriptscriptstyle Y_2$}  
        \put(20,20){$\scriptscriptstyle Y_r$}    
    \put(40,16){\circle*{3}}   \put(40,26){\circle*{3}}   \put(30,26){\circle*{3}}  \put(60,26){\circle*{3}}
       \put(40,16){\line(-1,1){10}}    \put(40,16){\line(0,1){10}}   \put(40,16){\line(2,1){20}}     
       \put(40,9){$\scriptscriptstyle a$}   \put(50,30){$\scriptscriptstyle \cdots$}   
       \put(26,30){$\scriptscriptstyle X_1$}   \put(38,30){$\scriptscriptstyle X_2$}  
        \put(58,30){$\scriptscriptstyle X_k$}     \put(10,6){\line(3,1){30}} 
      \end{picture}} \,,
\end{equation}
where $a$ is the largest label adjacent to the root such that $a<b$. Such an $a$ exists since $MD(T)$ is of size $deg(T)>1$. It is clear that $b$ is the largest value of the labels of $MD(T)$ (it is a decreasing tree, the root has the largest value). By choice, $a$ is the second largest value of the labels of $MD(T)$. We now consider the two subtrees   
$$
Z_1=\quad \raise -10pt\hbox{ \begin{picture}(40,30)
       \put(10,6){\circle{3}}   \put(0,16){\circle*{3}}   \put(-10,16){\circle*{3}}  \put(30,16){\circle*{3}}
       \put(10,6){\line(-2,1){20}}    \put(10,6){\line(-1,1){10}}   \put(10,6){\line(2,1){20}}     
       \put(10,-1){$\scriptscriptstyle b$}   \put(10,20){$\scriptscriptstyle \cdots$}   
       \put(-14,20){$\scriptscriptstyle Y_1$}   \put(-2,20){$\scriptscriptstyle Y_2$}  
        \put(28,20){$\scriptscriptstyle Y_r$}    
      \end{picture}} 
      \qquad\hbox{and}\qquad
Z_2= \quad \raise -10pt\hbox{ \begin{picture}(40,30)
       \put(10,6){\circle{3}}   \put(0,16){\circle*{3}}   \put(-10,16){\circle*{3}}  \put(30,16){\circle*{3}}
       \put(10,6){\line(-2,1){20}}    \put(10,6){\line(-1,1){10}}   \put(10,6){\line(2,1){20}}     
       \put(10,-1){$\scriptscriptstyle a$}   \put(10,20){$\scriptscriptstyle \cdots$}   
       \put(-14,20){$\scriptscriptstyle X_1$}   \put(-2,20){$\scriptscriptstyle X_2$}  
        \put(28,20){$\scriptscriptstyle X_k$}    
      \end{picture}} \,.
$$
Clearly $\deg(Z_1)<\deg(T)$ and $\deg(Z_2)<\deg(T)$. Hence by induction hypothesis we can find
a set partition $\Phi=\{\phi_1,\ldots,\phi_\ell\}\parti S$, a permutation $\sigma\colon[\ell]\to[\ell]$ and a family of trees  $\{T_i:T_i\in\calF[\phi_i]\}$ such that $\ T_{\sigma(1)}$ has root labelled by $b$ and
$T_{\sigma(k)}$ has root labelled by $a$ for some $k>1$. Furthermore $MD(Z_1)$ is the subtree of $Z_1$ labelled by the labels of the roots of $T_{\sigma(1)},\ldots,T_{\sigma(k-1)}$ and $MD(Z_2)$ is the subtree of $Z_2$ labelled by the labels of the roots of $T_{\sigma(k)},\ldots,T_{\sigma(\ell)}$. We can find this data in such a way that $Z_1$ is the leading term of $\Xi(\sb[T_{\sigma(1)}\cdots T_{\sigma(k-1)}])$ and $Z_2$ is the leading term of $\Xi(\sb[T_{\sigma(k)}\cdots T_{\sigma(\ell)}])$.

Using the same argument as before, it is clear that $T$ is the leading term of $[Z_1,Z_2]$. Thus, $T$ is the leading term of 
\begin{equation}\label{eq:sbT}
\big[\Xi(\sb[T_{\sigma(1)}\cdots T_{\sigma(k-1)}])\,,\,\Xi(\sb[T_{\sigma(k)}\cdots T_{\sigma(\ell)}])\big].
\end{equation}
We now need to show that the element in (\ref{eq:sbT}) is one of the element in the basis (\ref{eq:basis}).
This follows from the fact that $a$ is the second largest elements among the labels of the roots of $T_1,\ldots,T_\ell$. In particular it implies that the first step in the Lyndon bracketing of $\Xi(\sb[T_{\sigma(1)}\cdots T_{\sigma(\ell)}])$ is precisely the element in (\ref{eq:sbT}). Finally, $MD(T)$ is none other than the subtree of $T$ labelled by the labels of the roots of $T_{\sigma(1)},\ldots,T_{\sigma(\ell)}$. 
\end{proof}

\begin{exa} Let us compare the basis of $\calT[\{1,2,3\}]$ as given in
  Example~\ref{e:T3} with the following basis of 
$(\Lie\circ\calF)[\{1,2,3\}]$ as given by Eq.~(\ref{eq:basis}):
$$
 \raise -10pt\hbox{ \begin{picture}(30,30)
       \put(10,6){\circle{3}}   \put(0,16){\circle*{3}}   \put(20,16){\circle*{3}}
       \put(10,6){\line(-1,1){10}}    \put(10,6){\line(1,1){10}}
       \put(10,-1){$\scriptscriptstyle 1$}   \put(-4,20){$\scriptscriptstyle 2$}   \put(20,20){$\scriptscriptstyle 3$}  
      \end{picture}},\ 
 \raise -10pt\hbox{ \begin{picture}(12,30)
       \put(0,6){\circle{3}} \put(0,16){\circle*{3}} \put(0,26){\circle*{3}}
       \put(0,6){\line(0,1){20}}
       \put(3,5){$\scriptscriptstyle 1$} \put(3,17){$\scriptscriptstyle 2$} \put(3,27){$\scriptscriptstyle 3$}  
      \end{picture}} ,\ 
 \raise -10pt\hbox{ \begin{picture}(12,30)
       \put(0,6){\circle{3}}  \put(0,16){\circle*{3}} \put(0,26){\circle*{3}}
       \put(0,6){\line(0,1){20}}
       \put(3,5){$\scriptscriptstyle 1$}  \put(3,17){$\scriptscriptstyle 3$}  \put(3,27){$\scriptscriptstyle 2$}  
      \end{picture}} ,\ 
 \raise -10pt\hbox{ \begin{picture}(12,30)
       \put(0,6){\circle{3}}  \put(0,16){\circle*{3}} \put(0,26){\circle*{3}}
       \put(0,6){\line(0,1){20}} 
        \put(3,5){$\scriptscriptstyle 2$}  \put(3,17){$\scriptscriptstyle 3$} \put(3,27){$\scriptscriptstyle 1$}  
      \end{picture}} ,\ 
 \big[   \raise -7pt\hbox{ \begin{picture}(5,25)
       \put(0,6){\circle{3}} \put(0,16){\circle*{3}}  \put(0,6){\line(0,1){10}}
       \put(3,5){$\scriptscriptstyle 2$} \put(3,17){$\scriptscriptstyle 3$}   \end{picture}} \, ,\, 
         \raise -2pt\hbox{ \begin{picture}(5,5) \put(0,6){\circle{3}}  \put(3,5){$\scriptscriptstyle 1$}    
         \end{picture}}\,\,\big],\ 
 \big[      \raise -2pt\hbox{ \begin{picture}(5,5) \put(0,6){\circle{3}}  \put(3,5){$\scriptscriptstyle 2$}    
         \end{picture}}
          \, ,\,  \raise -7pt\hbox{ \begin{picture}(5,25)
       \put(0,6){\circle{3}} \put(0,16){\circle*{3}}  \put(0,6){\line(0,1){10}}
       \put(3,5){$\scriptscriptstyle 1$} \put(3,17){$\scriptscriptstyle 3$}   \end{picture}} \,\,\big],\ 
 \big[      \raise -2pt\hbox{ \begin{picture}(5,5) \put(0,6){\circle{3}}  \put(3,5){$\scriptscriptstyle 3$}    
         \end{picture}}
          \, ,\,  \raise -7pt\hbox{ \begin{picture}(5,25)
       \put(0,6){\circle{3}} \put(0,16){\circle*{3}}  \put(0,6){\line(0,1){10}}
       \put(3,5){$\scriptscriptstyle 1$} \put(3,17){$\scriptscriptstyle 2$}   \end{picture}} \,\,\big],\ 
 \big[      \raise -2pt\hbox{ \begin{picture}(5,5) \put(0,6){\circle{3}}  \put(3,5){$\scriptscriptstyle 3$}    
         \end{picture}}
          \, ,\,[ \raise -2pt\hbox{ \begin{picture}(5,5) \put(0,6){\circle{3}}  \put(3,5){$\scriptscriptstyle 2$}    
         \end{picture}} , \raise -2pt\hbox{ \begin{picture}(5,5) \put(0,6){\circle{3}}  \put(3,5){$\scriptscriptstyle 1$}    
         \end{picture}} \,\,]\big],\ 
 \big[  [    \raise -2pt\hbox{ \begin{picture}(5,5) \put(0,6){\circle{3}}  \put(3,5){$\scriptscriptstyle 3$}    
         \end{picture}}
        , \raise -2pt\hbox{ \begin{picture}(5,5) \put(0,6){\circle{3}}  \put(3,5){$\scriptscriptstyle 1$}    
         \end{picture}} \,\,   ] ,
         \raise -2pt\hbox{ \begin{picture}(5,5) \put(0,6){\circle{3}}  \put(3,5){$\scriptscriptstyle 2$}    
         \end{picture}} \,\,\big].
 $$
 The first four elements are already trees and they correspond to the basis of $\calF$ as given in Example~\ref{e:F3}.
As we expand the remaining elements in the basis  of trees (via $\Xi$) we get
    $$
     \big[   \raise -7pt\hbox{ \begin{picture}(5,25)
       \put(0,6){\circle{3}} \put(0,16){\circle*{3}}  \put(0,6){\line(0,1){10}}
       \put(3,5){$\scriptscriptstyle 2$} \put(3,17){$\scriptscriptstyle 3$}   \end{picture}} \, ,\, 
         \raise -2pt\hbox{ \begin{picture}(5,5) \put(0,6){\circle{3}}  \put(3,5){$\scriptscriptstyle 1$}    
         \end{picture}}\,\,\big]=
         \raise -10pt\hbox{  \begin{picture}(30,25)
       \put(10,6){\circle{3}} \put(0,16){\circle*{3}} \put(20,16){\circle*{3}}
       \put(10,6){\line(-1,1){10}}  \put(10,6){\line(1,1){10}}
       \put(10,-1){$\scriptscriptstyle 2$}   \put(-4,20){$\scriptscriptstyle 1$} \put(20,20){$\scriptscriptstyle 3$}  
      \end{picture}} +\ 
       \raise -10pt\hbox{ \begin{picture}(15,30)
       \put(0,6){\circle{3}}  \put(0,16){\circle*{3}} \put(0,26){\circle*{3}}
       \put(0,6){\line(0,1){20}} 
        \put(3,5){$\scriptscriptstyle 2$}  \put(3,17){$\scriptscriptstyle 3$} \put(3,27){$\scriptscriptstyle 1$}  
      \end{picture}} - \ 
       \raise -10pt\hbox{ \begin{picture}(15,30)
       \put(0,6){\circle{3}} \put(0,16){\circle*{3}} \put(0,26){\circle*{3}}
       \put(0,6){\line(0,1){20}}
       \put(3,5){$\scriptscriptstyle 1$} \put(3,17){$\scriptscriptstyle 2$} \put(3,27){$\scriptscriptstyle 3$}  
      \end{picture}},\qquad
   \big[      \raise -2pt\hbox{ \begin{picture}(5,5) \put(0,6){\circle{3}}  \put(3,5){$\scriptscriptstyle 2$}    
         \end{picture}}
          \, ,\,  \raise -7pt\hbox{ \begin{picture}(5,25)
       \put(0,6){\circle{3}} \put(0,16){\circle*{3}}  \put(0,6){\line(0,1){10}}
       \put(3,5){$\scriptscriptstyle 1$} \put(3,17){$\scriptscriptstyle 3$}   \end{picture}} \,\,\big]=
        \raise -10pt\hbox{ \begin{picture}(15,30)
       \put(0,6){\circle{3}}  \put(0,16){\circle*{3}} \put(0,26){\circle*{3}}
       \put(0,6){\line(0,1){20}}
       \put(3,5){$\scriptscriptstyle 2$}  \put(3,17){$\scriptscriptstyle 1$}  \put(3,27){$\scriptscriptstyle 3$}  
      \end{picture}} -  \raise -10pt\hbox{ \begin{picture}(30,30)
       \put(10,6){\circle{3}}   \put(0,16){\circle*{3}}   \put(20,16){\circle*{3}}
       \put(10,6){\line(-1,1){10}}    \put(10,6){\line(1,1){10}}
       \put(10,-1){$\scriptscriptstyle 1$}   \put(-4,20){$\scriptscriptstyle 2$}   \put(20,20){$\scriptscriptstyle 3$}  
      \end{picture}} 
-\  \raise -10pt\hbox{ \begin{picture}(15,30)
       \put(0,6){\circle{3}}  \put(0,16){\circle*{3}} \put(0,26){\circle*{3}}
       \put(0,6){\line(0,1){20}}
       \put(3,5){$\scriptscriptstyle 1$}  \put(3,17){$\scriptscriptstyle 3$}  \put(3,27){$\scriptscriptstyle 2$}  
      \end{picture}} ,
      $$
      $$
 \big[      \raise -2pt\hbox{ \begin{picture}(5,5) \put(0,6){\circle{3}}  \put(3,5){$\scriptscriptstyle 3$}    
         \end{picture}}
          \, ,\,  \raise -7pt\hbox{ \begin{picture}(5,25)
       \put(0,6){\circle{3}} \put(0,16){\circle*{3}}  \put(0,6){\line(0,1){10}}
       \put(3,5){$\scriptscriptstyle 1$} \put(3,17){$\scriptscriptstyle 2$}   \end{picture}} \,\,\big]
=       \raise -10pt\hbox{ \begin{picture}(15,30)
       \put(0,6){\circle{3}} \put(0,16){\circle*{3}}  \put(0,26){\circle*{3}}
       \put(0,6){\line(0,1){20}}
       \put(3,5){$\scriptscriptstyle 3$}  \put(3,17){$\scriptscriptstyle 1$}  \put(3,27){$\scriptscriptstyle 2$}  
      \end{picture}} 
- \   \raise -10pt\hbox{ \begin{picture}(30,30)
       \put(10,6){\circle{3}}   \put(0,16){\circle*{3}}   \put(20,16){\circle*{3}}
       \put(10,6){\line(-1,1){10}}    \put(10,6){\line(1,1){10}}
       \put(10,-1){$\scriptscriptstyle 1$}   \put(-4,20){$\scriptscriptstyle 2$}   \put(20,20){$\scriptscriptstyle 3$}  
      \end{picture}} 
- \raise -10pt\hbox{ \begin{picture}(15,30)
       \put(0,6){\circle{3}} \put(0,16){\circle*{3}} \put(0,26){\circle*{3}}
       \put(0,6){\line(0,1){20}}
       \put(3,5){$\scriptscriptstyle 1$} \put(3,17){$\scriptscriptstyle 2$} \put(3,27){$\scriptscriptstyle 3$}  
      \end{picture}} ,
      $$
      $$
   \big[      \raise -2pt\hbox{ \begin{picture}(5,5) \put(0,6){\circle{3}}  \put(3,5){$\scriptscriptstyle 3$}    
         \end{picture}}
          \, ,\,[ \raise -2pt\hbox{ \begin{picture}(5,5) \put(0,6){\circle{3}}  \put(3,5){$\scriptscriptstyle 2$}    
         \end{picture}} , \raise -2pt\hbox{ \begin{picture}(5,5) \put(0,6){\circle{3}}  \put(3,5){$\scriptscriptstyle 1$}    
         \end{picture}} \,\,]\big] 
         =      \raise -10pt\hbox{ \begin{picture}(15,30)
       \put(0,6){\circle{3}}  \put(0,16){\circle*{3}} \put(0,26){\circle*{3}}
       \put(0,6){\line(0,1){20}}
       \put(3,5){$\scriptscriptstyle 3$}   \put(3,17){$\scriptscriptstyle 2$}   \put(3,27){$\scriptscriptstyle 1$}  
      \end{picture}} 
 -  \raise -10pt\hbox{  \begin{picture}(30,30)
       \put(10,6){\circle{3}} \put(0,16){\circle*{3}} \put(20,16){\circle*{3}}
       \put(10,6){\line(-1,1){10}}  \put(10,6){\line(1,1){10}}
       \put(10,-1){$\scriptscriptstyle 2$}   \put(-4,20){$\scriptscriptstyle 1$} \put(20,20){$\scriptscriptstyle 3$}  
      \end{picture}} 
-  \raise -10pt\hbox{ \begin{picture}(15,30)
       \put(0,6){\circle{3}}  \put(0,16){\circle*{3}} \put(0,26){\circle*{3}}
       \put(0,6){\line(0,1){20}}
       \put(3,5){$\scriptscriptstyle 2$}  \put(3,17){$\scriptscriptstyle 1$}  \put(3,27){$\scriptscriptstyle 3$}  
      \end{picture}} 
-  \raise -10pt\hbox{ \begin{picture}(15,30)
       \put(0,6){\circle{3}} \put(0,16){\circle*{3}}  \put(0,26){\circle*{3}}
       \put(0,6){\line(0,1){20}}
       \put(3,5){$\scriptscriptstyle 3$}  \put(3,17){$\scriptscriptstyle 1$}  \put(3,27){$\scriptscriptstyle 2$}  
      \end{picture}} 
+   \raise -10pt\hbox{ \begin{picture}(30,30)
       \put(10,6){\circle{3}}   \put(0,16){\circle*{3}}   \put(20,16){\circle*{3}}
       \put(10,6){\line(-1,1){10}}    \put(10,6){\line(1,1){10}}
       \put(10,-1){$\scriptscriptstyle 1$}   \put(-4,20){$\scriptscriptstyle 2$}   \put(20,20){$\scriptscriptstyle 3$}  
      \end{picture}} 
+  \raise -10pt\hbox{ \begin{picture}(15,30)
       \put(0,6){\circle{3}} \put(0,16){\circle*{3}} \put(0,26){\circle*{3}}
       \put(0,6){\line(0,1){20}}
       \put(3,5){$\scriptscriptstyle 1$} \put(3,17){$\scriptscriptstyle 2$} \put(3,27){$\scriptscriptstyle 3$}  
      \end{picture}} ,
$$
      $$
       \big[  [    \raise -2pt\hbox{ \begin{picture}(5,5) \put(0,6){\circle{3}}  \put(3,5){$\scriptscriptstyle 3$}    
         \end{picture}}
        , \raise -2pt\hbox{ \begin{picture}(5,5) \put(0,6){\circle{3}}  \put(3,5){$\scriptscriptstyle 1$}    
         \end{picture}} \,\,   ] ,
         \raise -2pt\hbox{ \begin{picture}(5,5) \put(0,6){\circle{3}}  \put(3,5){$\scriptscriptstyle 2$}    
         \end{picture}} \,\,\big]
=    \raise -10pt\hbox{ \begin{picture}(30,30)
       \put(10,6){\circle{3}} \put(0,16){\circle*{3}}  \put(20,16){\circle*{3}}
       \put(10,6){\line(-1,1){10}}  \put(10,6){\line(1,1){10}}
       \put(10,-1){$\scriptscriptstyle 3$}   \put(-4,20){$\scriptscriptstyle 1$} \put(20,20){$\scriptscriptstyle 2$}  
      \end{picture}} 
+     \raise -10pt\hbox{ \begin{picture}(15,30)
       \put(0,6){\circle{3}} \put(0,16){\circle*{3}}  \put(0,26){\circle*{3}}
       \put(0,6){\line(0,1){20}}
       \put(3,5){$\scriptscriptstyle 3$}  \put(3,17){$\scriptscriptstyle 1$}  \put(3,27){$\scriptscriptstyle 2$}  
      \end{picture}} 
-  \raise -10pt\hbox{ \begin{picture}(15,30)
       \put(0,6){\circle{3}}  \put(0,16){\circle*{3}} \put(0,26){\circle*{3}}
       \put(0,6){\line(0,1){20}} 
        \put(3,5){$\scriptscriptstyle 2$}  \put(3,17){$\scriptscriptstyle 3$} \put(3,27){$\scriptscriptstyle 1$}  
      \end{picture}} 
-  \raise -10pt\hbox{ \begin{picture}(30,30)
       \put(10,6){\circle{3}}   \put(0,16){\circle*{3}}   \put(20,16){\circle*{3}}
       \put(10,6){\line(-1,1){10}}    \put(10,6){\line(1,1){10}}
       \put(10,-1){$\scriptscriptstyle 1$}   \put(-4,20){$\scriptscriptstyle 2$}   \put(20,20){$\scriptscriptstyle 3$}  
      \end{picture}} 
-  \raise -10pt\hbox{ \begin{picture}(15,30)
       \put(0,6){\circle{3}}  \put(0,16){\circle*{3}} \put(0,26){\circle*{3}}
       \put(0,6){\line(0,1){20}}
       \put(3,5){$\scriptscriptstyle 1$}  \put(3,17){$\scriptscriptstyle 3$}  \put(3,27){$\scriptscriptstyle 2$}  
      \end{picture}} 
+  \raise -10pt\hbox{ \begin{picture}(15,30)
       \put(0,6){\circle{3}}  \put(0,16){\circle*{3}} \put(0,26){\circle*{3}}
       \put(0,6){\line(0,1){20}}
       \put(3,5){$\scriptscriptstyle 2$}  \put(3,17){$\scriptscriptstyle 1$}  \put(3,27){$\scriptscriptstyle 3$}  
      \end{picture}} 
. $$
We then remark that each tree of the basis in Example~\ref{e:T3} appear once as the leading term (the first term) of an expression above.
\end{exa}

\begin{rems}
We have that the Lie bracket on $\calT$ is filtrated with respect to our degree. 
That is if we start with two disjoint sets $I,J$ (each with a linear order) and two elements in 
$T\in\calT[I]$ and $Y\in\calT[J]$ where $\deg(T)=d_1$ and $\deg(Y)=d_2$, then the maximal degree part 
of $[T,Y]$ is of degree $d_1+d_2$. This follows from the fact that $MD(Z)$ of the  term $Z$ appearing in the maximal degree part of $[T,Y]$ must be obtained from a subtree of either the grafting of $MD(T)$ in $MD(Y)$ or the other way around. Grafting at the root will produce a decreasing tree in one case (hence achieving the degree $d_1+d_2$). In general, it is possible that other grafting of $MD(T)$ in $MD(Y)$ (or the other way around) achieve the maximal degree. This did not happen in the proof of the Theorem~\ref{th:TLieF1} because the label of the roots of $MD(T)$ and $MD(Y)$ where the largest two labels.
\end{rems}

\section{$\calT = \Lie \circ \calF$ as species}\label{S:LieFsp}

In the previous section we have shown that $\calT[S] = \Lie \circ \calF[S]$ as vector space for any finite set $S$ with a linear order. Using the basis~(\ref{eq:basis}) in Theorem~\ref{th:TLieF1} we can now define the notion of Lie-degree. We say that an element $\Xi(\sb[T_{\sigma(1)}\cdots T_{\sigma(\ell)}])$ has Lie-degree $\ell=\hbox{L-deg}(\Xi(\sb[T_{\sigma(1)}\cdots T_{\sigma(\ell)}]))$. For an arbitrary element $\varphi\in\calT[S]$ we let $\hbox{L-deg}(\varphi)=\ell$ if $\ell$ is the smallest Lie-degree of the basis elements with non-zero coefficient in the expansion of $\varphi$ in the basis~\ref{eq:basis}.
It is clear that given two disjoint sets $I,J$ with linear order on each, 
$\varphi\in\calT[I]$ and $\psi\in\calT[J]$, we have that 
$\hbox{L-deg}([\varphi,\psi])=\hbox{L-deg}(\varphi)+\hbox{L-deg}(\psi)$ for any linear order on $I+J$ compatible with the orders on $I$ and $J$. The notion of Lie-degree is not related to the notion of degree we used in the proof of Theorem~\ref{th:TLieF1}.

 We consider the specie $[\calT,\calT]$ which is the image of $[,\,]\colon\calT\bullet\calT\to\calT$. 
It follows from our discussion above that $\calT[S]\Big/ \lower2pt\hbox{$[\calT,\calT][S]$}$ will be isomophic to the set of homogeneous elements of Lie-degree equal to one. These elements are precisely $\calF[S]$. We thus have the following corollary:

\begin{coro}\label{co:FT} There is a linear isomorphism
  $\calF[S] \to \calT[S]\Big/ \lower2pt\hbox{$[\calT,\calT][S]$}$.
\end{coro}

%
%

Since $\calT\Big/ \lower2pt\hbox{$[\calT,\calT]$}$ is clearly a specie, Corollary~\ref{co:FT} allows us to view $\calF$ as a specie.
 \begin{exa} \label{ex:F2}
 Using the basis~(\ref{eq:basis}) of $\calT[3]$, the space $\calF[3]=\calT[3]\Big/ \lower2pt\hbox{$[\calT,\calT][3]$}$ is the linear span of the elements:
   $$
 \raise -10pt\hbox{ \begin{picture}(30,30)
       \put(10,6){\circle{3}}   \put(0,16){\circle*{3}}   \put(20,16){\circle*{3}}
       \put(10,6){\line(-1,1){10}}    \put(10,6){\line(1,1){10}}
       \put(10,-1){$\scriptscriptstyle 1$}   \put(-4,20){$\scriptscriptstyle 2$}   \put(20,20){$\scriptscriptstyle 3$}  
      \end{picture}},\ 
 \raise -10pt\hbox{ \begin{picture}(12,30)
       \put(0,6){\circle{3}} \put(0,16){\circle*{3}} \put(0,26){\circle*{3}}
       \put(0,6){\line(0,1){20}}
       \put(3,5){$\scriptscriptstyle 1$} \put(3,17){$\scriptscriptstyle 2$} \put(3,27){$\scriptscriptstyle 3$}  
      \end{picture}} ,\ 
 \raise -10pt\hbox{ \begin{picture}(12,30)
       \put(0,6){\circle{3}}  \put(0,16){\circle*{3}} \put(0,26){\circle*{3}}
       \put(0,6){\line(0,1){20}}
       \put(3,5){$\scriptscriptstyle 1$}  \put(3,17){$\scriptscriptstyle 3$}  \put(3,27){$\scriptscriptstyle 2$}  
      \end{picture}} ,\ 
 \raise -10pt\hbox{ \begin{picture}(12,30)
       \put(0,6){\circle{3}}  \put(0,16){\circle*{3}} \put(0,26){\circle*{3}}
       \put(0,6){\line(0,1){20}} 
        \put(3,5){$\scriptscriptstyle 2$}  \put(3,17){$\scriptscriptstyle 3$} \put(3,27){$\scriptscriptstyle 1$}  
      \end{picture}} ,\ 
 \big[   \raise -7pt\hbox{ \begin{picture}(5,25)
       \put(0,6){\circle{3}} \put(0,16){\circle*{3}}  \put(0,6){\line(0,1){10}}
       \put(3,5){$\scriptscriptstyle 2$} \put(3,17){$\scriptscriptstyle 3$}   \end{picture}} \, ,\, 
         \raise -2pt\hbox{ \begin{picture}(5,5) \put(0,6){\circle{3}}  \put(3,5){$\scriptscriptstyle 1$}    
         \end{picture}}\,\,\big],\ 
 \big[      \raise -2pt\hbox{ \begin{picture}(5,5) \put(0,6){\circle{3}}  \put(3,5){$\scriptscriptstyle 2$}    
         \end{picture}}
          \, ,\,  \raise -7pt\hbox{ \begin{picture}(5,25)
       \put(0,6){\circle{3}} \put(0,16){\circle*{3}}  \put(0,6){\line(0,1){10}}
       \put(3,5){$\scriptscriptstyle 1$} \put(3,17){$\scriptscriptstyle 3$}   \end{picture}} \,\,\big],\ 
 \big[      \raise -2pt\hbox{ \begin{picture}(5,5) \put(0,6){\circle{3}}  \put(3,5){$\scriptscriptstyle 3$}    
         \end{picture}}
          \, ,\,  \raise -7pt\hbox{ \begin{picture}(5,25)
       \put(0,6){\circle{3}} \put(0,16){\circle*{3}}  \put(0,6){\line(0,1){10}}
       \put(3,5){$\scriptscriptstyle 1$} \put(3,17){$\scriptscriptstyle 2$}   \end{picture}} \,\,\big],\ 
 \big[      \raise -2pt\hbox{ \begin{picture}(5,5) \put(0,6){\circle{3}}  \put(3,5){$\scriptscriptstyle 3$}    
         \end{picture}}
          \, ,\,[ \raise -2pt\hbox{ \begin{picture}(5,5) \put(0,6){\circle{3}}  \put(3,5){$\scriptscriptstyle 2$}    
         \end{picture}} , \raise -2pt\hbox{ \begin{picture}(5,5) \put(0,6){\circle{3}}  \put(3,5){$\scriptscriptstyle 1$}    
         \end{picture}} \,\,]\big],\ 
 \big[  [    \raise -2pt\hbox{ \begin{picture}(5,5) \put(0,6){\circle{3}}  \put(3,5){$\scriptscriptstyle 3$}    
         \end{picture}}
        , \raise -2pt\hbox{ \begin{picture}(5,5) \put(0,6){\circle{3}}  \put(3,5){$\scriptscriptstyle 1$}    
         \end{picture}} \,\,   ] ,
         \raise -2pt\hbox{ \begin{picture}(5,5) \put(0,6){\circle{3}}  \put(3,5){$\scriptscriptstyle 2$}    
         \end{picture}} \,\,\big].
 $$
 The first four elements form a basis of $\calF[3]$ and the remaining ones are zero modulo $[\calT,\calT]$.
 The action of the symmetric group $S_3$ is given by the action on the quotient. For example if we have the transposition $\sigma=(1\, 2)$ then 
 $$\sigma\big(\,  \raise -10pt\hbox{ \begin{picture}(25,30)
       \put(10,6){\circle{3}}   \put(0,16){\circle*{3}}   \put(20,16){\circle*{3}}
       \put(10,6){\line(-1,1){10}}    \put(10,6){\line(1,1){10}}
       \put(10,-1){$\scriptscriptstyle 1$}   \put(-4,20){$\scriptscriptstyle 2$}   \put(20,20){$\scriptscriptstyle 3$}  
      \end{picture}} \big) = 
       \raise -10pt\hbox{ \begin{picture}(30,30)
       \put(10,6){\circle{3}}   \put(0,16){\circle*{3}}   \put(20,16){\circle*{3}}
       \put(10,6){\line(-1,1){10}}    \put(10,6){\line(1,1){10}}
       \put(10,-1){$\scriptscriptstyle 2$}   \put(-4,20){$\scriptscriptstyle 1$}   \put(20,20){$\scriptscriptstyle 3$}  
      \end{picture}} =  \big[   \raise -7pt\hbox{ \begin{picture}(5,25)
       \put(0,6){\circle{3}} \put(0,16){\circle*{3}}  \put(0,6){\line(0,1){10}}
       \put(3,5){$\scriptscriptstyle 2$} \put(3,17){$\scriptscriptstyle 3$}   \end{picture}} \, ,\, 
         \raise -2pt\hbox{ \begin{picture}(5,5) \put(0,6){\circle{3}}  \put(3,5){$\scriptscriptstyle 1$}    
         \end{picture}}\,\,\big]\  -\   \raise -10pt\hbox{ \begin{picture}(12,30)
       \put(0,6){\circle{3}}  \put(0,16){\circle*{3}} \put(0,26){\circle*{3}}
       \put(0,6){\line(0,1){20}} 
        \put(3,5){$\scriptscriptstyle 2$}  \put(3,17){$\scriptscriptstyle 3$} \put(3,27){$\scriptscriptstyle 1$}  
      \end{picture}} 
+\   \raise -10pt\hbox{ \begin{picture}(12,30)
       \put(0,6){\circle{3}}  \put(0,16){\circle*{3}} \put(0,26){\circle*{3}}
       \put(0,6){\line(0,1){20}} 
        \put(3,5){$\scriptscriptstyle 1$}  \put(3,17){$\scriptscriptstyle 2$} \put(3,27){$\scriptscriptstyle 3$}  
      \end{picture}}       
      \equiv\ 
-\   \raise -10pt\hbox{ \begin{picture}(12,30)
       \put(0,6){\circle{3}}  \put(0,16){\circle*{3}} \put(0,26){\circle*{3}}
       \put(0,6){\line(0,1){20}} 
        \put(3,5){$\scriptscriptstyle 2$}  \put(3,17){$\scriptscriptstyle 3$} \put(3,27){$\scriptscriptstyle 1$}  
      \end{picture}} 
+\   \raise -10pt\hbox{ \begin{picture}(12,30)
       \put(0,6){\circle{3}}  \put(0,16){\circle*{3}} \put(0,26){\circle*{3}}
       \put(0,6){\line(0,1){20}} 
        \put(3,5){$\scriptscriptstyle 1$}  \put(3,17){$\scriptscriptstyle 2$} \put(3,27){$\scriptscriptstyle 3$}  
      \end{picture}}       .
 $$
 \end{exa}

Once we have the identification of $\calF=\calT\Big/ \lower2pt\hbox{$[\calT,\calT]$}$, we then have a natural action of the symmetric group on $\calF[n]$ and on $\Lie\circ\calF[n]$. To see the action of $S_n$ on $\Lie\circ\calF[n]$, we consider the natural order on $[n]$ and the basis
\begin{equation}\label{eq:basisLF}
      \big\{\sb[T_{\sigma(1)}\cdots T_{\sigma(\ell)}]:\substack{\Phi=\{\phi_1,\ldots,\phi_\ell\}\parti S,\
               \sigma\colon[\ell]\to[\ell],\ T_i\in\calF[\phi_i],\ T_{\sigma(1)} \hbox{ \small has the largest root}}
      \big\}\,.      \end{equation}
A permutation $\pi\in S_n$ acts on a basis element $\sb[T_{\sigma(1)}\cdots T_{\sigma(\ell)}]$ as follow.  First by acting on each trees $\pi(T_{\sigma(1)}),\ldots ,\pi(T_{\sigma(\ell)})$. We have to rewrite the 
$\pi(T_{\sigma(i)})$ as linear combination of trees in
$\calF[\pi(\phi_i)]$ where $\phi_i$ are the set of labels of
$T_{\sigma(i)}$.  We then substitute the results in
$\sb[T_{\sigma(1)}\cdots T_{\sigma(\ell)}]$. We use the Jacobi
relation and 
antisymmetry to rewrite the result as a linear combination of elements in the 
basis~(\ref{eq:basisLF}). The important fact to notice is that the element $\pi\big(\sb[T_{\sigma(1)}\cdots T_{\sigma(\ell)}]\big)$ will be a linear combination of elements of the basis~(\ref{eq:basisLF}) with exactly the same Lie-degree (the same number of $\calF$-trees bracketed). Hence the matrix representation corresponding to $\pi\in S_n$ acting on the basis~(\ref{eq:basisLF}) is block diagonal, each block corresponds to the Lie-degrees. 

On the other hand, we have shown in the  Theorem~\ref{th:TLieF1} that under the map  $\Xi\colon\Lie\circ \calF[S]\to\calT[S]$ the basis~(\ref{eq:basisLF}) is map to the basis
\begin{equation}\label{eq:basis2}
\big\{\Xi(\sb[T_{\sigma(1)}\cdots T_{\sigma(\ell)}]):\substack{\Phi=\{\phi_1,\ldots,\phi_\ell\}\parti S,\
               \sigma\colon[\ell]\to[\ell],\ T_i\in\calF[\phi_i],\ T_{\sigma(1)} \hbox{ \small has the largest root}}
      \big\}\,.
      \end{equation}
Now a permutation $\pi\in S_n$ acts on a basis element $\Xi(\sb[T_{\sigma(1)}\cdots T_{\sigma(\ell)}])$ as follows.  First by acting on each trees $\pi(T_{\sigma(1)}),\ldots ,\pi(T_{\sigma(\ell)})$. We have to rewrite the 
$\pi(T_{\sigma(i)})$ as linear combination of basis elements in $\calT[\pi(\phi_i)]$ where $\phi_i$ are the set of labels of $T_{\sigma(i)}$.  We remark that the Lie-degree of the result will also be one but will contains higher Lie-degree terms. We then substitute the results in $\sb[T_{\sigma(1)}\cdots T_{\sigma(\ell)}]$. We use the Jacobi relation and antisymmetry to rewrite the result as a linear combination of element in the 
basis~(\ref{eq:basis2}). The important fact to notice in this case is that the element $\pi\big(\Xi(\sb[T_{\sigma(1)}\cdots T_{\sigma(\ell)}])\big)$ will be a linear combination of elements of the basis~(\ref{eq:basisLF}) with Lie-degree equal or higher than $\Xi(\sb[T_{\sigma(1)}\cdots T_{\sigma(\ell)}])$. Hence the matrix representation corresponding to $\pi\in S_n$ acting on the basis~(\ref{eq:basis2}) is block triangular, each block correspond to the Lie-degrees. Moreover the block diagonal part of this matrix is exactly the same as the action of $\pi$ on $\Lie\circ\calF$. We have thus shown that the $S_n$ module $\calT[n]$ is $S_n$-isomorphic to $\Lie\circ\calF[n]$. This shows the following corollary:

\begin{coro} \label{co:TLieF2}
    $\calT[S] = \Lie\circ \calF[S]$ as species. 
\end{coro}

\section{Operations on $ \calF$ }\label{S:F}

We recall first the operad structure on $\calT$ as explained in \cite{CL01}.
To give the operadic structure it is enough to explain the composition
on two elements, that is the compositions
$$\circ_i: \calT[I]\otimes \calT[J]\rightarrow \calT[I\setminus\{i\}+J],$$
for two disjoint sets $I$ and $J$ and for $i\in I$.
Let $T\in\calT[I]$, $In(T,i)$ the set of incoming edges at the vertex
labelled by $i$ in $T$. The composition is defined by
$$T\circ_i S=\sum_{f:In(T,i)\rightarrow J} T\circ_i^f S$$
where $ T\circ_i^f S$ is the rooted tree obtained by substituting the tree $S$ for the vertex $i$ in $T$. 
The outgoing edge of $i$, if it exists, becomes the outgoing edge of the root of $S$, whereas the 
incoming edges of $i$ are graphting
on the vertices of $S$ following the map $f$. The root is the root of $T$ or the root of $S$ if $i$ is the 
root of $T$. Here is an example:


\bigskip

$$  \raise -10pt\hbox{ \begin{picture}(25,30)
       \put(0,16){\circle*{3}}    \put(20,16){\circle*{3}}     \put(10,6){\circle{3}}
       \put(0,16){\line(1,-1){10}}    \put(20,16){\line(-1,-1){10}}                                 
       \put(-3,20){$\scriptscriptstyle a$}  \put(22,20){$\scriptscriptstyle b$}  \put(10,-1){$\scriptscriptstyle i$} 
       \end{picture}}\ \circ_i  
\raise -10pt\hbox{ \begin{picture}(25,30)
       \put(10,16){\circle*{3}}     \put(10,6){\circle{3}}
       \put(10,16){\line(0,-1){10}}    
       \put(10,20){$\scriptscriptstyle \alpha$}  \put(10,-1){$\scriptscriptstyle \beta$} 
       \end{picture}} =
     \raise -10pt\hbox{ \begin{picture}(25,30)
       \put(10,16){\circle*{3}}   \put(0,26){\circle*{3}}   \put(20,26){\circle*{3}}      \put(10,6){\circle{3}}
       \put(10,16){\line(-1,1){10}}    \put(10,16){\line(1,1){10}}                           \put(10,6){\line(0,1){10}}        
       \put(12,13){$\scriptscriptstyle \alpha$}   \put(-2,30){$\scriptscriptstyle a$}  
        \put(20,30){$\scriptscriptstyle b$}                                           \put(10,-1){$\scriptscriptstyle \beta$} 
       \end{picture}}\, 
       +\ \  \raise -10pt\hbox{ \begin{picture}(15,30)
       \put(0,16){\circle*{3}}   \put(0,26){\circle*{3}}   \put(10,16){\circle*{3}}      \put(0,6){\circle{3}}
       \put(0,16){\line(0,1){10}}    \put(0,6){\line(1,1){10}}                           \put(0,6){\line(0,1){10}}        
       \put(-8,13){$\scriptscriptstyle \alpha$}   \put(-2,30){$\scriptscriptstyle a$}  
        \put(10,20){$\scriptscriptstyle b$}                                           \put(0,-1){$\scriptscriptstyle \beta$} 
       \end{picture}}\ 
 +\  \ \raise -10pt\hbox{ \begin{picture}(15,30)
       \put(0,16){\circle*{3}}   \put(0,26){\circle*{3}}   \put(10,16){\circle*{3}}      \put(0,6){\circle{3}}
       \put(0,16){\line(0,1){10}}    \put(0,6){\line(1,1){10}}                           \put(0,6){\line(0,1){10}}        
       \put(-8,13){$\scriptscriptstyle \alpha$}   \put(-2,30){$\scriptscriptstyle b$}  
        \put(10,20){$\scriptscriptstyle a$}                                           \put(0,-1){$\scriptscriptstyle \beta$} 
       \end{picture}}\
  + \raise -10pt\hbox{ \begin{picture}(25,30)
       \put(0,16){\circle*{3}}   \put(10,16){\circle*{3}}   \put(20,16){\circle*{3}}      \put(10,6){\circle{3}}
       \put(0,16){\line(1,-1){10}}    \put(20,16){\line(-1,-1){10}}    \put(10,6){\line(0,1){10}}        
       \put(8,20){$\scriptscriptstyle \alpha$}   \put(-3,20){$\scriptscriptstyle a$}  
        \put(21,20){$\scriptscriptstyle b$}                                           \put(10,-1){$\scriptscriptstyle \beta$} 
       \end{picture}}\,  $$ 

\bigskip

If $I$ and $J$ are endowed with a linear order then $I\setminus\{i\}+J$ is endowed with the order on $I$ and $J$
and for all $x$ in $I$ one has $x<J$ if and only if $x<i$ and $x>J$ if and only if $x>i$. Recall that the degree
of a tree $T$ is  the number of vertices of $MD(T)$, 
the maximal decreasing connected subtree of $T$ from the root. Given an ordered set $I$, the vector space
$\calT[I]$ is filtered by the degree: $F_d\calT[I]$ is spanned by the trees $T$ of order less than $d$.
Therefore $\calF[I]=F_1\calT[I]$.

\begin{theo}Let $I,J$ be two ordered sets and let $i\in I$. The composition
$$\circ_i: \calT[I]\otimes\calT[J]\rightarrow \calT[I\setminus\{i\}+J]$$
maps $F_d\calT[I]\otimes F_e\calT[J]$ to $F_{d+e-1}\calT[I\setminus\{i\}+J].$
\end{theo}

\begin{proof}
Let $T\in\calT[I]$, $S\in\calT[J]$ and $f:In(T,i)\rightarrow J$.

If $i$ is a vertex of $MD(T)$ then any vertex $x$ in $T$ lying in the path from the root of $T$ to $i$ 
satisfies $x>i$. Then $x>J$ and $x$ is a vertex of $MD(T\circ_i^f S)$. 
Also any vertex of $MD(S)$ is a vertex of $MD(T\circ_i^f S)$.
Let $y$ be a vertex of $T$ attached to $i$ by an edge 
$E$ in $In(T,i)$. If $y$ is a vertex of $MD(T)$ then $y<i$ and $y<J$. If $f$ sends $E$ to a vertex 
of $MD(S)$ then $y$ is a vertex of $MD(T\circ_i^f S)$. If it doesn't then
the degree of $T\circ_i^f S$ is strictly less than $d+e-1$. If $y$ is not a vertex of $MD(T)$ then 
$y>i$ and $y>J$, hence $y$ is not a vertex of  $MD(T\circ_i^f S)$.

 It is also clear that any other vertex of $T$ which is not in $MD(T)$ or any vertex which is
not in $MD(S)$ won't lie in $MD(T\circ_i^f S)$. As a consequence the degree of 
$T\circ_i^f S$ is less than $d+e-1$. It is exactly $d+e-1$ if $f$ sends any vertex of $MD(T)$ attached to
$i$ by an edge in $In(T,i)$ to a vertex of $MD(S)$. In this case $MD(T\circ_i^f S)=MD(T)\circ_i^{\tilde f} MD(S)$
where $\tilde f$ is the restriction of $f$ to the set of edges of $MD(T)$.

If $i$ is not a vertex of $MD(T)$ then one sees easily that $MD(T\circ_i^f S)=MD(T)$. 
\end{proof}

\begin{coro} The collection $(\calF[n])_{n\geq 1}$ forms a sub nonsymmetric operad of the pre-Lie operad $\calT$.
\end{coro}

\begin{proof}We apply the previous theorem to $d=e=1$.
\end{proof}


 \end{document}